\documentclass{amsart}
\usepackage{amsmath}
\usepackage{amssymb}
\usepackage{amsthm,amssymb,amsmath,enumerate,graphicx}
\usepackage{color,soul}
\usepackage{subfigure}
\usepackage{caption}

\makeatletter
\newtheorem*{rep@theorem}{\rep@title}
\newcommand{\newreptheorem}[2]{%
\newenvironment{rep#1}[1]{%
 \def\rep@title{#2 \ref{##1}}%
 \begin{rep@theorem}}%
 {\end{rep@theorem}}}
\makeatother

\newtheorem{theorem}{Theorem}[section]

\newreptheorem{theorem}{Theorem}

\newtheorem{lemma}[theorem]{Lemma}

\newtheorem{conjecture}[theorem]{Conjecture}

\theoremstyle{remark}
\newtheorem{notation}[theorem]{Notation}
\newtheorem{remark}[theorem]{Remark}

\newcommand{\cc}{\mathcal C}

\newcommand{\cf}{\mathcal F}
\newcommand{\cg}{\mathcal G}

\newcommand{\ci}{\mathcal I}

\begin{document}

\makeatletter

\makeatother
\author{Ron Aharoni}
\address{Department of Mathematics\\ Technion, Haifa, Israel}
\email[Ron Aharoni]{raharoni@gmail.com}
\thanks{\noindent The research of the first author was
supported by BSF grant no. $2006099$, by  ISF grant no. $1581/12$,
by the Technion's research promotion fund, and by
the Discont Bank chair.}

\author{Eli Berger}
\address{Department of Mathematics, Haifa University, Haifa 31999, Israel}
\email[Eli Berger]{berger.haifa@gmail.com}
\thanks{\noindent The research of the second author was
supported by BSF grant no. $2006099$, and by  ISF grant no. $1581/12$.}

\author{Dani Kotlar}
\address{Computer Science Department, Tel-Hai college, Upper Galilee, Israel}
\email[Dani Kotlar]{dannykotlar@gmail.com}

\author{Ran Ziv}
\address{Computer Science Department, Tel-Hai college, Upper Galilee, Israel}
\email[Ran Ziv]{ranzivziv@gmail.com}

\title{Degree conditions for matchability in $3$-partite hypergraphs}
\maketitle
\begin{abstract}
We study  conjectures relating degree conditions in $3$-partite hypergraphs to the matching number of the hypergraph, and
 use topological methods to prove special cases. In particular,
we prove a strong version of a theorem of Drisko
\cite{drisko} (as generalized by the first two authors \cite{ab}), that every
family of $2n-1$ matchings of size $n$ in a bipartite graph has a
partial rainbow matching of size $n$. We show that  milder restrictions on the sizes of the matchings suffice.
Another
result that is strengthened is a theorem of Cameron and Wanless \cite{CamWan}, that every
Latin square has a diagonal (permutation submatrix) in which no
symbol appears more than twice. We show that the same is true under
the weaker condition that the square is row-Latin.
\end{abstract}


\section{Rainbow matchings and matchings in $3$-partite $3$-uniform hypergraphs}

All  conjectures and results mentioned in this paper can be  traced back to a by now well known conjecture of Ryser \cite{ryser}, that for $n$ odd every
Latin square possesses a full transversal, namely a permuatation submatrix with distinct symbols.   Brualdi \cite{brualdi} and Stein \cite{stein} conjectured that for general $n$,
every $n \times n$ Latin square possesses
a partial transversal of size $n-1$. Stein \cite{stein} generalized this still further,  replacing the Latinity condition by the milder requirement that each of the
$n$ symbols appears in $n$ cells  (implying, among other things, that each cell contains
precisely one symbol). See \cite{wansurvey} for a survey on these conjectures and  related results.

These conjectures can be formulated in the terminology of $3$-partite hypergraphs.  We say that a pair $(X,Y)$ of sides of  a $k$-partite hypergraph
$H$ is {\em simple} if no pair $(x,y)$ with $x \in X, y \in Y$
appears in more than one edge of $H$. A hypergraph is called {\em
simple} if no edge repeats more than once.
In this terminology, an $n \times n$  Latin
square is an $n$-regular $3$-partite hypergraph with all sides  of
size $n$, and all three pairs of sides being simple. For a hypergraph $H$ denote by $\nu(H)$  the maximal size of a matching
in $H$.  In this language the Brualdi-Stein conjecture is:

\begin{conjecture}\label{conj:rbs}
Let $H$ be an $n$-regular $3$-partite hypergraph with sides
$A,B, C$, all of size $n$, and assume that all three pairs of sides,
$(A,B)$, $(B,C)$ and $(A,C)$, are simple. Then $\nu(H) \ge n-1$.
\end{conjecture}

While Stein's conjecture is:

\begin{conjecture}\label{conj:stein}
Let $H$ be an $n$-regular $3$-partite hypergraph with sides
$A,B, C$, all of size $n$, and assume that the pair
$(A,B)$  is simple. Then $\nu(H) \ge n-1$.
\end{conjecture}

In fact, we believe that an even stronger (and more simply stated) conjecture is true:

\begin{conjecture}\label{symmetricstein}
Let $H$ be a simple $n$-regular $3$-partite hypergraph with sides
of size $n$. Then $\nu(H) \ge n-1$.
\end{conjecture}

In the final section of this paper we shall present some yet more general conjectures. But let us first turn
to a concept pertinent to these conjectures, that of {\em rainbow matchings}.
Given a family $F_1, \ldots, F_m$ of sets, a choice $f_1 \in F_1,
\ldots ,f_m \in F_m$ is called a {\em rainbow set}. If the elements
of the sets $F_i$ are themselves sets, and if the sets $f_i$ are
disjoint, then the rainbow set is  called a {\em rainbow matching}.
If the elements of the sets in $F_i$ are edges of a bipartite graph $G$, then a
rainbow matching is a matching in a $3$-partite hypergraph, in which
the vertices of one side represent the sets $F_i$, and the other two
sides are those of $G$.

A generalization of Conjecture \ref{conj:stein} was proposed in \cite{ab}:

\begin{conjecture}\label{conj:ab}
Any $n$ matchings of size $n$ in a bipartite graph have a partial rainbow matching of size $n-1$.
\end{conjecture}

Many partial results have been obtained on this conjecture, see, e.g.,  \cite{woolbright, shorhatami, ach, kotlarziv, kotlarzivaharoni, germans, cep, Pok}.

What happens if we demand that $\nu \ge n$, rather than $n-1$, or in the terminology of rainbow matchings we want a rainbow matching of size $n$? Strangely,
we need to almost double the requirement on the number of matchings. Drisko \cite{drisko} proved a theorem which was later
generalized in \cite{ab} to:

\begin{theorem}\label{thm:drisko}
A family of $2n-1$ matchings of size $n$ in a bipartite graph has a
rainbow matching of size $n$.
\end{theorem}

This is sharp - repeating $n-1$ times each of the matchings consisting respectively of the even edges in $C_{2n}$ and the odd edges in $C_{2n}$ shows that $2n-2$ matchings do not suffice. In \cite{akzequality} this was shown to be the only example demonstrating the sharpness of the theorem.

  In $3$-partite
hypergraphs terminology, Theorem \ref{thm:drisko} reads:

\begin{theorem}
Let $H$ be a $3$-partite hypergraph with sides $A,B, C$, and assume
that
\begin{enumerate}
\item $|A|=2n-1$.
\item $deg(a)=n$ for every $a \in A$.
\item The pairs  $(A,B)$ and $(A,C)$ are simple.
\end{enumerate}
Then $\nu(H) \ge n$.
\end{theorem}

It is plausible that condition (3) in this theorem is too
restrictive.

\begin{conjecture}\label{conj:drisko}
Let $H$ be a  $3$-partite hypergraph with sides $A,B, C$, and
assume that
\begin{enumerate}
\item $|A|\ge 2n-1$.
\item $deg(a)\ge n$ for every $a \in A$.
\item $deg(v) \le 2n-1$ for every $v \in B \cup C$.
\end{enumerate}
Then $\nu(H) \ge n$.
\end{conjecture}

The proofs of Theorem \ref{thm:drisko} given in \cite{drisko, ab} were combinatorial. A topological proof
 yields a stronger version:

\begin{theorem}\label{thm:drisko_improved}
Let $\cg=(V,E)$ be a bipartite graph and let
$\cf=\{F_1,F_2,\ldots,F_{2n-1}\}$ be a family of matchings in $\cg$
so that $|F_i|\ge i$ for $i=1,\ldots,n-1$ and $|F_i|=n$ for
$i=n,\ldots,2n-1$. Then $\cf$ has a rainbow matching of size $n$.
\end{theorem}

In fact,  this condition is also necessary. Call a
sequence of numbers $a_1 \le a_2\le \ldots \le a_{2n-1}$ {\em
accommodating} if  every family $\cf=\{F_1,F_2,\ldots,F_{2n-1}\}$
 of matchings in a bipartite graph satisfying $|F_i|\ge a_i$ has a
 rainbow matching of size $n$.

 \begin{theorem}\label{thm:accommodating}
A sequence $a_1 \le a_2\le \ldots \le a_{2n-1}$ is accommodating if
and only if $a_i \ge \min(i, n)$  for all $i \le 2n-1$.
 \end{theorem}

We call a pair $(X,Y)$ of sides in a $k$-partite hypergraph $H$ \emph{$p$-simple} if no pair $(x,y)$, for $x \in X, y \in Y$ is contained in more than $p$ edges of $H$.

The following is  a special case of Conjecture \ref{conj:drisko}:

\begin{theorem}\label{almostconjdrisko}
Let $H$ be a  $3$-partite hypergraph with sides $A,B, C$, and
assume that
\begin{enumerate}
\item $|A|\ge 2n-1$ and $|B|=|C|=n$,
\item $deg(a)= n$ for every $a \in A$, and
\item The pair $(A,C)$ is simple and the pair $(B,C)$ is 2-simple.
\end{enumerate}
Then $\nu(H) \ge n$.\end{theorem}

This theorem will yield a strengthening of the following result of
Cameron and Wanless \cite{CamWan}, which can be regarded as ``half''  of Conjecture
\ref{conj:stein}:

\begin{theorem}\label{thm:camwan:original}
Every Latin square contains a permutation submatrix in which no
symbol appears more than twice.
\end{theorem}

 For a set $S$ of vertices in a hypergraph $H$ let $\delta(S)$ (resp.
$\Delta(S)$) be the minimal (resp. maximal) degree of a vertex in
$S$. In this terminology,
Theorem \ref{thm:camwan:original} is:

\begin{theorem}\label{cawan3partite}
Let $H$ be an $n$-regular $3$-partite hypergraph with  sides $A,B,C$
such that $|A|=|B|=|C|=n$, and all pairs of sides are simple. Then
there exists a set $F$ of $n$ edges, satisfying
\begin{enumerate}
\item

$\Delta_F(B \cup C)=\delta_F(B \cup C)=1$
\item
$\Delta_F(A) \le 2$.
\end{enumerate}
\end{theorem}

Here the vertices in $A$ represent symbols, those in $B$ represent
columns, and those in $C$ represent rows. Theorem
\ref{cawan3partite} follows quite directly from Theorem
\ref{thm:drisko}:

\begin{proof}
 For each $a \in A$ let $M_a=N_a =\{(b,c) \mid b \in B, c \in C \mid
(a,b,c) \in H\}$. By the simplicity assumption, these are matchings of size $n$. By Theorem~\ref{thm:drisko},
the set of matchings $\{M_a \mid a \in A\} \cup \{N_a \mid a \in A\}$
possesses a rainbow matching $R$ of size $n$. Defining $F=\{(a,b,c) \mid (b,c) \in
R\}$ yields the desired result.
\end{proof}

The topological tools will allow us to strengthen Theorem \ref{thm:camwan:original} to row-Latin
squares, namely squares in which no symbol appears twice in the same row (but may appear more than once in the same column).
In $3$-partite hypergraph terminology:

\begin{theorem}\label{strongcawan}
For the conclusion of Theorem \ref{cawan3partite} to hold it
suffices that the two pairs $(A,C)$ and $(B,C)$ are simple.
\end{theorem}

Of course, by symmetry, it suffices also to assume that $(A,B)$ and
$(B,C)$ are simple.

\section{A topological tool}\label{topology}

For a graph $G$ denote by $\ci(G)$ the complex (closed down
hypergraph) of  independent sets in $G$. If $G=L(H)$, the line graph
of a hypergraph $H$, then $\ci(G)$ is the complex of matchings in
$H$.
 A simplicial complex  $\cc$ is called (homologically)
\emph{$k$-connected}  if for every $-1 \le j \le k$, the $j$-th reduced simplicial homology group of $\cc$ with rational coefficients $\tilde{H}_j(\cc)$ vanishes. The (homological) \emph{connectivity} $\eta_H(\cc)$ is the largest $k$ for which $\cc$ is $k$-connected, plus $2$.

\begin{remark}
\noindent
\begin{itemize}
\item[(a)] This is a shifted (by $2$) version of the usual definition of connectivity. The shift simplifies the statements below, as well as the statements of basic properties of the connectivity parameter.
\item[(b)] If $\tilde{H}_j(\cc)=0$ for all $j$ then we define  $\eta_H(\cc)=\infty$.
\item[(c)] There exists also a homotopical notion of connectivity, $\eta_h(\cc)$: it is the minimal dimension of a ``hole'' in the complex.  The first topological version of Hall's theorem \cite{ah} used that notion.
The relationship between the two parameters is that $\eta_H\ge \eta_h$ for all complexes, and if $\eta_h(\cc)\ge 3$, meaning that the complex is simply connected, then $\eta_H(\cc)= \eta_h(\cc)$.  All facts mentioned in this article (in particular, the main tool we use, the Meshulam game) apply also to $\eta_h$.
\end{itemize}
\end{remark}

\begin{notation}
  Given sets $(V_i)_{i=1}^n$  and a subset  $I$ of $[n]=\{1,\ldots,n\}$, we write $V_I$ for $\bigcup_{i\in
    I}V_i$. For $A \subseteq V(G)$ we denote by $\ci(G) \upharpoonright A$ the complex of independent sets in the graph induced by $G$ on $A$.
\end{notation}

Given  sets $(V_i)_{i=1}^n$ of  vertices in a graph $G$, an {\em independent transversal}   is an
independent set containing at least one vertex from each $V_i$. Note
that in this definition the transversal needs not be injective,
namely a set $V_i$ may be represented twice, which makes a difference
if the sets are not disjoint. In our application (as is the case in
the  most common applications of the theorem) the sets $V_i$ are
disjoint, in which case the transversal may well be assumed to be
injective.

 The following is a topological version of Hall's theorem:

\begin{theorem}\label{ah}
  If $\eta_H(\ci(G) \upharpoonright V_I) \ge |I|$ for every $I \subseteq
  [n]$ then there exists an independent transversal.
\end{theorem}

Variants of this theorem appeared implicitly in \cite{ah} and
\cite{me1},  and the theorem is stated and proved explicitly as
Proposition~1.6 in \cite{me2}.

A standard argument of adding dummy vertices yields the  deficiency
version of Theorem \ref{ah}:
\begin{theorem}\label{deficiency}
If $\eta_H(\ci(G[\bigcup_{i \in I}V_i])) \ge |I|-d$ for every $I
\subseteq [m]$ then the system has a partial  independent
transversal of size $m-d$.
\end{theorem}

In order to apply these theorems, combinatorially formulated lower
bounds on $\eta_H(\ci(G))$ are needed. One such bound is  due to
Meshulam \cite{me2} and is conveniently expressed in terms of a game
between two players, CON and NON, on the graph $G$. CON wants to
show high connectivity, NON wants to thwart her attempt. At each
step, CON chooses an edge $e$ from the graph remaining at this
stage, where in the first step the graph is $G$. NON can then either

\begin{enumerate}
\item
Delete $e$ from the graph. We call such a step a ``deletion'', and denote the resulting graph by ``$G-e$''.

or

\item Remove the two endpoints of $e$, together
  with all neighbors of these vertices and the edges incident to them, from the graph. We call such a step
  an ``explosion'', and denote the resulting graph by ``$G*e$''.

\end{enumerate}

\noindent The result of the game (payoff to CON) is defined as follows: if at some point there
remains an isolated vertex, the result is $\infty$. Otherwise, at some
point all vertices have disappeared, in which case the result of the
game is the number of explosion steps. We define $\Psi(G)$ as the value of the game, i.e., the maximum payoff to CON in an optimal play of both players.

\begin{theorem}\label{etaPsi} \cite{me2}
$\eta_H(\ci(G))\ge \Psi(G)$.
\end{theorem}

\begin{remark}
 This formulation of $\Psi$ is equivalent to a recursive definition of $\Psi(G)$ as the maximum over all edges of $G$, of  $\min(\Psi(G-e),\Psi(G*e)+1)$.  For an explicit proof of Theorem~\ref{etaPsi} using the recursive definition of $\Psi$, see Theorem~1 in \cite{adba}. The ``game''  formulation first appeared in \cite{abz}.
\end{remark}

We shall use the Meshulam bound on line graphs. If $G=L(H)$ then playing the game means that CON
offers NON a pair of edges in $H$ having a common vertex. Deletion in $G$
corresponds to separating the two edges at the vertex where they
meet. Explosion corresponds to removing the three endpoints of these
edges.

\section{Proof of Theorem \ref{thm:accommodating}}\label{drisko}

\begin{lemma}\label{lemma:1}
Let $G$ be a bipartite graph with sides $U,W$, and assume that there exist
 $u_1, \ldots ,u_{2\ell-1} \in U$ satisfying
 $deg(u_i)\ge \min (i,\ell)$ for all $i \le 2\ell-1$ . Then $\Psi(L(G))\ge \ell$.
\end{lemma}

\begin{proof}
By induction on $\ell$. For $\ell=1$ the statement is obvious. For
$\ell>1$ CON chooses an edge $e=u_1w$ for some $w \in W$, and offers
NON in any order all pairs $e,f$ for edges $f=yw\neq e$ (here $y \in
U$). Assume first that $NON$ explodes one of these pairs.
 Since $G$ is simple, the degree of each vertex in
$U\cap V(G')$ is reduced by the explosion by at most $1$. This
implies that the remaining graph $G'$ satisfies the hypothesis of
the theorem with  $\ell-1$ replacing $\ell$. By the induction
hypothesis, we have $\Psi(L(G'))\ge \ell-1$ and hence $\Psi(L(G))\ge
\ell$.

Next consider the case that NON separates $e$ from $f$ for all
$f=yw\neq e$. Then CON offers all pairs $e,g$ for  $g=u_1z \neq e$
(here $z \in W$). NON has to explode one of these pairs, say $e,g$
for $g=u_1z$,  so as not to render $e$ isolated. Since all pairs
$e,f$ for $f=yw\neq e$ have been separated, this explosion preserves
the vertex $w$, removing at it only the edge $u_1w$. Thus, again,
the degrees of all $u_i, i>1$ are reduced by only $1$, and the
condition of the lemma holds with $\ell-1$ replacing $\ell$. Again,
the desired conclusion follows by induction.
\end{proof}

\begin{proof}[Proof of Theorem~\ref{thm:accommodating}]
For the sufficiency part, let $a_1 \le a_2\le \ldots \le a_{2n-1}$
be  a sequence of natural numbers satisfying $a_i \ge \min(i, n)$
for all $i \le 2n-1$ and let  $\cf=\{F_1,F_2,\ldots,F_{2n-1}\}$ be a
family
 of matchings in a bipartite graph $\cg=(V,E)$ satisfying $|F_i|\ge a_i$ for all $i \le 2n-1$.
 We need to show that $\cf$ has a rainbow matching of size $n$.
Let $A$ and $B$ be the two sides of $V$ and let $m=|A|\ge n$.  By
considering a third side $C$ of size $2n-1$ whose vertices
correspond to the matchings $F_i$ we obtain a 3-partite 3-hypergraph
$H$. We need to show that $H$ has a matching of size $n$.

Consider the bipartite graph $G$ induced on $B$ and $C$.  Since the
sets $F_i$ are matchings,  $G$ is simple. The vertices in $A$ induce
a partition $V_1, \ldots ,V_m$ on $E(G)$. By the hypothesis of the
theorem, the vertices of $C$ can be ordered as $c_i, ~i \le 2n-1$,
where $deg_G(c_i) \ge \min(i,n)$. For $I \subseteq A$ let $G^I$ be
the graph with sides $B,C$ induced by $I$, namely having as edges
pairs $(b,c)$ completed by some $a \in I$ to an edge of $H$. Write
$k=|A|-|I|$. Since the sets $F_i$ are matchings, meaning that the
pair of sides $(A,C)$ in $H$ is simple, for every $c \in C$ we have
$deg_{G^I}(c) \ge deg_G(c)-k$. This entails that the conditions of
Lemma \ref{lemma:1} are valid for $G^I$ with $\ell=n-k$. By the
lemma, we have $\eta_H(\ci(L(G^I))) \ge n-k$.
 By Theorem \ref{deficiency} this  suffices to
show that $\nu(H) \ge n$, namely there is a rainbow matching of size
$n$.

For the other direction of the theorem, let $a_1,\ldots,a_n$  be an
ascending sequence of natural numbers such that  $a_k\le k-1$ for some $k\le n$. Let $C_{2n}$ be a  cycle of size
$2n$ and let $M\cup N$ be a partition of its edges into two
matchings, each of size $n$. Consider the following family
$\cf=\{F_1,F_2,\ldots,F_{2n-1}\}$ of matchings. Each of
$F_{n+1},\ldots,F_{2n-1}$ is a  copy of $M$. Let
$N=\{e_1,\ldots,e_n\}$ and for each $i=1,\ldots,n$ let $F_i$ be a
copy of the matching $\{e_1,\dots,e_{\min(n,a_i)}\}$. Clearly, these
matchings satisfy the condition of the theorem. We claim that they
do not possess a rainbow matching of size $n$. Clearly, a matching
of size $n$ is contained either in $M$ or in $N$, and since there
are only $n-1$ matchings $F_i$ that meet $M$, a rainbow matching of
size $n$ must represent the matchings $F_1, \ldots ,F_n$. But this
is impossible, since the union of the matchings $F_1, \ldots ,F_k$
contain together fewer than $k$ edges.
\end{proof}

It may be worth noting that an anaolgous version of Conjecture \ref{conj:ab} fails. Let $Q^1, \ldots ,Q^{2k}$ be disjoint copies of $P_3$, the path with three
edges, let $O_i$ be the set of two odd  edges in $Q^i~~(i \le 2k)$,  let $e_i$ be the middle edge in $Q^i$, let $F_i=\{e_1, \ldots ,e_k\}$ for $i \le k$, and let
$F_i=\bigcup\{O_j \mid 1 \le j \le i-k\} \cup \{e_j \mid j>i-k\}$ for $i>k$. Then $|F_i| \ge i$ for all $i$, and the largest rainbow matching is of size $\lfloor \frac{3k}{2}\rfloor$.
\section{Proof of Theorems \ref{almostconjdrisko} and \ref{strongcawan}}\label{CamWan}

For the convenience of the reader, let us repeat Theorem \ref{almostconjdrisko}:

\begin{reptheorem}{almostconjdrisko}

Let $H$ be a  $3$-partite hypergraph with sides $A,B, C$, and
assume that
\begin{enumerate}
\item $|A|\ge 2n-1$ and $|B|=|C|=n$,
\item $deg(a)= n$ for every $a \in A$, and
\item The pair $(A,C)$ is simple and the pair $(B,C)$ is 2-simple.
\end{enumerate}
Then $\nu(H) \ge n$.
\end{reptheorem}


\begin{proof}
Since the pair $(A,C)$ is simple, it follows from  (2)  that the degree of each vertex in $C$ is $|A|$.
Since $(B,C)$ is 2-simple, we have
\begin{equation} \label{deltas}
\Delta(B)\le 2n.
\end{equation}

For  $b \in B$ let
$V_b=\{(a,c) \mid (a,b,c) \in E(H)\}$ be a set of edges in $A \times C$. We need to show that the
sets $V_b$ have a full rainbow matching.
For $K\subseteq B$  let $G_K$ be the  bipartite graph with sides $A$ and $C$, and with edge set $\bigcup_{b \in K}V_b$.
By Theorems \ref{ah} and \ref{etaPsi} it suffices to show that

\begin{equation} \label{etak}
\Psi(\ci(L(G_K))) \ge |K|
\end{equation}
for every set $K \subseteq B$. Write $k=|K|$. The minimal value of $|E(G_K)|$ occurs when all the vertices in $B\setminus K$ have maximal degree, which
by (\ref{deltas}) means:

\begin{equation}\label{minEK}
|E(G_K)|\ge n(2n-1)-2n(n-k)=2nk-n,
\end{equation}

and by the 2-simplicity of the pair $(B,C)$ we have

\begin{equation}\label{DeltaGK}
\Delta_{G_K}(C)\le 2k.
\end{equation}

We play Meshulam's game on $G_K$ as follows. Let $u_1,u_2,\ldots, u_n$ be the vertices of $C$. For each $i \le n$ let $d_i$ be the degree of $u_i$ in $G_K$ and assume that $d_1\le d_2\le\cdots\le d_n$. CON goes over the vertices in this order. Let $\ell_i$ be the degree of $u_i$ at the time it is handled. For each $i$, if $\ell_i\ge 2$ then CON offers all pairs of edges meeting at $u_i$, in any order. This he does  until NON explodes a pair, or until all  edges meeting at $u_i$ are separated from each other. If $\ell_i< 2$ then $u_i$ is skipped and CON handles the next vertex in the list.

Let $p_i$ be the number of explosions performed until $u_i$ was handled, including the possible explosion at $u_i$ itself. Assume first that for some $i$ NON separated all pairs of edges meeting at $u_i$ and $p_i+\ell_i\ge k$. Let $e_1,\dots,e_{\ell_i}$ be the edges meeting at $u_i$ at the time it is handled and let $w_1,\ldots,w_{\ell_i}$ be their corresponding endpoints at $A$. For each  $j=1,\ldots,\ell_i$, CON offers the pairs of edges $(e_j,f)$ for all edges  $f$ meeting $e_j$ at $w_j$. Note that at least one such $f$ exists for each $w_j$, otherwise $e_j$ is isolated, meaning that the score of the game is $\infty$. Also note that NON must explode a pair in each $w_j$, otherwise the corresponding $e_j$ will become isolated. Thus, CON scores $\ell_i$ points at $u_i$, which together with the $p_i$ already scored, the score of the game is at least $k$.

Hence, if we make the negation assumption that the score of the game is less than $k$, then for each $u_i$ one of the following two occurs:

\begin{enumerate}
  \item [POS1:] NON exploded a pair at $u_i$, or
  \item [POS2:] NON separated all the edges meeting at $u_i$ and $p_i+\ell_i< k$.
\end{enumerate}

In each explosion two vertices from $A$ are removed along with their incident edges. Thus, as $G_K$ is simple, in each explosion the degree of each vertex in $C$ decreases by at most two. Hence,

\begin{equation}\label{ell_i:lowerbound}
    \ell_i\ge d_i-2p_i\textrm { for all } i.
\end{equation}

Suppose NON separated all the edges meeting at $u_i$. Then, by POS2 we have $p_i+\ell_i< k$. This together with (\ref{ell_i:lowerbound}) yield,

\textbf{Claim 1:} If NON separated all the edges meeting at $u_i$ then $d_i < k+p_i$.

For each $i=1,\ldots,n$ let $\pi_i=p_i+(n-i)$. Note that by the negation assumption $\pi_n=p_n<k$, so there exists  a  minimal index $t$ for which $\pi_t<k$.

\textbf{Claim 2:} NON separated all the edges meeting at $u_t$.

\begin{proof}[Proof of Claim 2]
Consider first the case $t=1$. If NON exploded a pair at $u_1$ then $p_1=1$ and we have $\pi_1=n\ge k$, contradicting the definition of $t$. So we may assume that $t>1$. If NON exploded a pair at $u_t$ then $p_{t-1}=p_t-1$ and thus $\pi_{t-1}=p_{t-1}+(n-(t-1))=p_t-(n-t)<k$, contradicting the minimality of $t$.
\end{proof}

Now, by the minimality of $t$ we have $\pi_{t-1}=p_{t-1}+(n-(t-1))\ge k$ and $\pi_t=p_t+(n-t)< k$. By Claim 2, $p_{t-1}=p_t$. Hence $p_{t-1}+(n-(t-1))= k$ and $p_t+(n-t)= k-1$.
Thus

\textbf{Claim 3:} $t=n-k+p_t+1$.

We calculate an upper bound on $|E(G_K)|$. By Claims 1 and 2 and the fact that the $d_i$s are  ascending  we have
$d_i<k+p_t$ for all $i=1,\ldots,t$. From this and Claim 3, we conclude that the first $t$ vertices are incident to less than $(n-k+p_t+1)(k+p_t)$ edges in $|E(G_K)|$. By (\ref{DeltaGK}) and Claim 3, the remaining $n-t$ edges are incident to at most $(k-p_t-1)2k$ edges. So, we have,

\begin{equation}\label{EGK:upperbound}
    |E(G_K)| < (n-k+p_t+1)(k+p_t)+(k-p_t-1)2k.
\end{equation}

Let $s=k-p_t$. Then (\ref{EGK:upperbound}) can be written in a somewhat simpler form:

\begin{equation}\label{EGK:upperbound:simpler}
    |E(G_K)| < (n-s+1)(2k-s)+(s-1)2k.
\end{equation}

Let $m=s(n-s+1)$. Rearranging terms in (\ref{EGK:upperbound:simpler}) we obtain

\begin{equation}\label{EKUpperBound}
    |E(G_K)| < 2nk-m,
\end{equation}

By the negation assumption $s\ge 1$ implying $m\ge n$. By (\ref{EKUpperBound}) it follows that $|E(G_K)| < 2nk-n$, contradicting (\ref{minEK}).
\end{proof}

\begin{proof}[Proof of Theorem~\ref{strongcawan}]
Let $A'$ be the union of two identical copies of $A$, that is, $A'=A \cup A^\dagger$, where $A^\dagger=\{a^\dagger \mid a \in A\}$ and let $H'$ be the hypergraph with sides $A',B,C$, defined by $E(H')=E(H) \cup \{(a^\dagger,b,c) \mid (a,b,c) \in H\}$.
We have $|A'|=2n$. Also, since $(A,C)$ is simple so is $(A',C)$, and since $(B,C)$ is simple, the pair $(B,C)$ is 2-simple in $H'$. By Theorem~\ref{almostconjdrisko}, we have $\nu(H')=n$, which implies the desired result.
\end{proof}

\section{Possible generalizations}

In \cite{kotlarzivaharoni} the following conjecture was proposed:

\begin{conjecture}\label{fracd}
Let $H$ be  a simple  $3$-partite $d$-regular  hypergraph  with sides
of size $n$.

\begin{enumerate}
\item
If $d \le n$ then $\nu(H) \ge \frac{d-1}{d}n$.
 \item
 If $d \ge 2n-1$ then $\nu(H) =n$.
 \end{enumerate}
\end{conjecture}

Part (1) would imply Conjecture \ref{symmetricstein}. Part (2) is sharp. To see this,
let $a,b,c$ be vertices in the respective sides $A,B,C$ of a hypergraph $H$ with $|A|=|B|=|C|=n$, put in $E(H)$ the set
$\{(a,b,x) \mid x \in C\setminus \{c\}\} \cup \{(a,y,c) \mid y \in B\setminus \{b\}\} \cup \{(z,b,c) \mid z \in A\setminus \{a\}\}$, and complete it to
a $2n-2$ regular hypergraph by adding edges not containing any of $a,b,c$.
In such a hypergraph $\nu \le n-1$, since $a,b,c$ cannot be covered by the same matching.

  An asymmetric formulation of the conjecture may better capture its essence:

\begin{conjecture}
Let  $H$ be  a simple  $3$-partite hypergraph
with sides $A,B,C$.

\begin{enumerate}
\item
If $d=\delta(A ) \ge \Delta(B \cup C)$ then $\nu(H) \ge
\frac{d-1}{d}|A|$.

\item
If $\delta(A ) \ge \max(\Delta(B \cup C), 2|A|-1)$ then
$\nu(H)=|A|$.
\end{enumerate}
\end{conjecture}

\begin{remark} Theorem \ref{ah} can be used to prove that if $\delta(A) \ge 2\Delta(B
\cup C)-1$ then $\nu(H)=|A|$.
\end{remark}

 Note that in item (2)
 there is a jump by a factor of $2$ with respect to (1), similar to that between Conjecture \ref{conj:ab} and Theorem \ref{thm:drisko}. By (1) to get $\nu(H)\ge |A|-1$ we only (conjecturally) need
$\delta(A ) \ge \max(\Delta(B \cup C), |A|)$.

A conjecture generalizing Theorem \ref{thm:drisko} in the same
spirit is:

\begin{conjecture}
Let  $H$ be  a simple  $3$-partite hypergraph
with sides $A,B,C$, and suppose that $|A| =2n-1$, $deg(a)\ge n$ for
all $a \in A$, and $deg(v) \le 2n-1$ for all $v \in B \cup C$. Then
$\nu(H) \ge n$.
\end{conjecture}


\end{document}